\documentclass[a4paper,11pt]{article}
\usepackage{amssymb,palatino}
\newtheorem{theorem}{Theorem}

\newtheorem{proposition}{Proposition}
\newcommand{\sty}{\displaystyle}
\newcommand{\QED}{\begin{flushright} $\Box$ \end{flushright}}
\newcommand{\binomiale}[2]{\left( \begin{array}{c} #1\\#2 \end{array} \right)}

\newcommand{\legendre}[2]{( #1 | #2 )}
\begin{document}

\title{Additive decompositions induced \\ by multiplicative characters over finite fields}
\author{Davide Schipani\thanks{University of Zurich, Switzerland}, ~~
    Michele Elia\thanks{Politecnico di Torino, Italy}
}

\maketitle

\thispagestyle{empty}

\begin{abstract}
\noindent

In 1952, Perron showed that quadratic residues in a field of prime order
satisfy certain additive properties. This result has been generalized in different directions, and 
our contribution is to provide a further generalization concerning multiplicative quadratic and cubic characters over any finite field. 
In particular, recalling that a character partitions the multiplicative group of the field into cosets with respect 
 to its kernel, we will derive 
 the number of representations of an element as a sum of two elements belonging to two given cosets. These numbers are then related to the equations satisfied by the polynomial characteristic functions of the cosets. 

Further, we show a connection, a quasi-duality, with the problem of determining how many elements can be added to each element of a subset of a coset in such a way as to obtain elements still belonging to a subset of a coset. 

\end{abstract}



\paragraph{Keywords:} Characters, Residuacity, Finite Fields

\vspace{2mm}
\noindent
{\bf Mathematics Subject Classification (2010): } 12Y05, 12E20, 12E30 

\vspace{8mm}

\section{Introduction}

Back in 1952, Perron \cite{perron} proved that every quadratic residue in $\mathbb F_p$, with $p$ an odd prime, can be written as the sum of two quadratic residues in exactly $\lfloor\frac{p+1}{4}\rfloor-1$ ways, and as the sum of two quadratic non-residues in exactly $\lfloor\frac{p+1}{4}\rfloor$ ways, a symmetric statement also holding for quadratic non-residues.

Winterhof \cite{winter} generalized this result, proving that if $x_j$ is a nonzero element in a finite field $\mathbb F_q$, $\chi$ a nontrivial multiplicative character of order $n$ and $\omega$ a primitive $n$-th root of unity, then 
$\sigma_0+1=\sigma_1=\cdots=\sigma_{n-1}=\frac{q-1}{n}$, where $\sigma_i$, $i=0,\ldots,n-1$, stands for the number of field elements $x$ such that $\chi(x)\bar\chi(x+x_j)=\omega^i$.

A generalization in a different direction was provided by Monico and Elia \cite{elia2, monico}: they proved that the only partition of the field satisfying the additive properties found by Perron is that of quadratic residues and non-residues, providing a sort of converse and a way to define residues additively; moreover, they generalized Perron's result for any multiplicative character defined over a prime field.

Lastly, in his dissertation \cite{raymond} Raymond showed how to generalize the case of the quadratic character to any finite field.

We provide here further generalizations considering quadratic and cubic characters over any finite field. Section 2 provides notations and preliminaries for the rest of the paper. In Section 3 we derive 
expressions for the number of (ordered) representations of an element as a sum of two elements belonging to two given cosets of the character partition. Furthermore, these expressions are shown to be related to the equations satisfied by the polynomial characteristic functions of the cosets. 
In Section 4 we point out a sort of duality relationship with another problem, and mention how the two problems are involved in polynomial factorization. 

\section{Preliminaries}


Let $\mathbb F_{p^m}$, $p$ an odd prime, be a finite field with polynomial basis 
 $\{ 1,\eta, \eta^2, \ldots , \eta^{m-1} \}$ where $\eta$ is a root of an irreducible polynomial of degree
 $m$ over $\mathbb F_p$, and let $\mathcal B_0$ be the set of squares (excluding $0$) and $\mathcal B_1$
 the complementary set in $\mathbb F_{p^m}^*$. The quadratic character is a mapping from 
 $\mathbb F_{p^m}^*$ into the complex numbers defined as
$$  \chi_2(\alpha^h \theta)  = (-1)^h ~~~~\theta\in\mathcal B_0,~~~~h=0,1~~,  $$ 
 where $\alpha$ is a primitive element in $\mathbb F_{p^m}^*$. Furthermore, we set $\chi_2(0)=0$. 

We can define an indicator function of the sets $\mathcal B_j$ using the
  quadratic character, namely, for every $\gamma \neq 0$,
%
$$    I_{\mathcal B_j}(\gamma) = \frac{1+(-1)^{j}\chi_2(\gamma)}{2} = \left\{
     \begin{array}{l}
       1 ~~~~\mbox{if}~~ \gamma \in \mathcal B_j   \\
       0 ~~~~\mbox{otherwise}
     \end{array} \right.
        ~~j=0,1~~. $$ 
Also, writing $\gamma=\gamma_0+\gamma_1 \eta + \cdots + \gamma_{m-1}\eta^{m-1}$, a polynomial characteristic function that identifies the set $\mathcal B_j$ can be defined as
$$  f_{\mathcal B_j}(\mathbf X) = \sum_{\gamma \in \mathbb F_{q}^*} I_{\mathcal B_j}(\gamma) \mathbf X^{\gamma}  ~~, $$
where $\mathbf X^{\gamma}$ is a short notation standing for the monomial 
 $ x_0^{\gamma_0}x_1^{\gamma_1}\cdots x_{m-1}^{\gamma_{m-1}} $. Note that this notation allows us
  to formally write
$$ \mathbf X^{\gamma} \mathbf X^{\delta}= \mathbf X^{\gamma+\delta} ~~~~ \bmod \mathfrak I_{\mathbf X} ~~,    $$
where $\mathfrak I_{\mathbf X}=\langle (x_1^p-1), \ldots , (x_m^p-1) \rangle$ is an ideal  
 in $\mathbb Q[\mathbf X]$.
  
It is immediate to see that, if $\Phi_p(x)$ is the $p$-th cyclotomic polynomial,
\begin{equation}
 \label{coeff0}
  f_{\mathcal B_0}(\mathbf X)+f_{\mathcal B_1}(\mathbf X)  +1 = \sum_{\gamma \in \mathbb F_{p^m}^*} \mathbf X^{\gamma}+1 = \prod_{i=0}^{m-1} \Phi_p(x_i) ~~. 
\end{equation}

\noindent
In the following we will indicate the last product with the notation $\Phi(\mathbf X)$.
\vspace{3mm}

In the case of cubic characters, non-trivial characters exist only if $p$ is $2$ with even 
 exponent $m$, or $p$ is an odd prime congruent to $1$ modulo $6$, or $p$ is an odd prime congruent to
 $5$ modulo $6$ with even exponent $m$. 

Let us write any nonzero element of the field as $\alpha^{k+3n}$, with $k\in\{0,1,2\}$ and $\alpha$ a primitive element: we define $\mathcal A_0=\{\alpha^{3i}:~~i=0,\ldots , \frac{p^m-1}{3}-1 \}$,
 that is the subgroup of the cubic powers in $\mathbb F_{p^m}^*$, and let $\mathcal A_1=\alpha \mathcal A_0$
 and $\mathcal A_2=\alpha^2 \mathcal A_0$ be the two cosets that complete the coset 
 partition of the set of nonzero elements of $\mathbb F_{p^m}$.

Similarly to the case of the quadratic character, we define an indicator function of the sets $\mathcal A_j$
 using a cubic character, that is a mapping of the type
$$  \chi_3(\alpha^h \theta) =  \omega^h ~~~~\theta\in\mathcal A_0,~~~ h=0,1,2~~,  $$ 
with $\omega$ a primitive cubic root of unity in $\mathbb C$ (and $\chi_3(0)=0$ by definition).

The indicator function, for every $x \neq 0$, is then
$$    I_{\mathcal A_j}(x) = \frac{1+\omega^{2j}\chi_3(x)+\omega^{j}\bar \chi_3(x)}{3} = \left\{
     \begin{array}{l}
       1 ~~~~\mbox{if}~~ x \in \mathcal A_j   \\
       0 ~~~~\mbox{otherwise}
     \end{array} \right.
        ~~j=0,1,2~~, $$ 
(where the bar denotes complex conjugation), and a characteristic function that identifies
 the set $\mathcal A_j$ can be defined in the same way as above:
$$  f_{\mathcal A_j}(\mathbf X) = \sum_{\gamma \in \mathbb F_{p^m}^*} I_{\mathcal A_j}(\gamma) \mathbf X^{\gamma}  ~~~~j=0,1,2 ~~. $$
Again it is immediate to see that 
\begin{equation}
 \label{coeff2}
  f_{\mathcal A_0}(\mathbf X)+f_{\mathcal A_1}(\mathbf X)+f_{\mathcal A_2}(\mathbf X)   +1 = \sum_{\gamma \in \mathbb F_{p^m}^*} \mathbf X^{\gamma}+1 = \Phi(\mathbf X) ~~. 
\end{equation}

\section{Results}
To count the number of representations of a $\beta\neq 0$ in the field $\mathbb F_{p^m}$
 as the sum of an element in $\mathcal B_j$ and an element in $\mathcal B_i$ ($i$ not necessarily
 different  from $j$), it suffices to compute
\begin{equation}
 \label{count2}
\sum_{z\neq 0,\beta}\frac{1+(-1)^{j}\chi_2(z)}{2}~~\frac{1+(-1)^{i}\chi_2(\beta -z)}{2}
\end{equation}

Analogously, when we have a cubic character, to count the number of representations of a $\beta\neq 0$ in the field $\mathbb F_{p^m}$
 as the sum of an element in $\mathcal A_j$ and an element in $\mathcal A_i$ ($i$ not necessarily different
 from $j$), it suffices to compute
\begin{equation}
 \label{count3}
\sum_{z\neq 0,\beta}\frac{1+\omega^{2j}\chi_3(z)+\omega^{j}\bar \chi_3(z)}{3}~~\frac{1+\omega^{2i}\chi_3(\beta-z)+\omega^{i}\bar \chi_3(\beta-z)}{3}.
\end{equation}
We summarize the conclusions in the next three theorems.

\begin{theorem}
  \label{lemma2}
The number of representations 
$R_{p^m}^{(2)}(\beta,i,j)$ of a $\beta\neq 0$ in the field $\mathbb F_{p^m}$,
 $p$ an odd prime, as the sum of an element in $\mathcal B_j$ and an element in $\mathcal B_i$ is
\begin{equation}
 \label{coeff1}
R_{p^m}^{(2)}(\beta,i,j) =\frac{1}{4}\left(p^m-2-\chi_2(\beta)(-1)^i-\chi_2(\beta)(-1)^j-(-1)^{i+j}\chi_2(-1)\right)~~,
\end{equation}
 and depends only on the quadratic residuacity of $\beta$.
\end{theorem}

\noindent
{ \sc Proof}. The proof is immediate from (\ref{count2}), since,
 for any nontrivial character $\chi$, $\sum_{x \in \mathbb F_{p^m}} \chi(x)=0$ and
 $\sum_{x \in \mathbb F_{p^m}} \chi(x) \bar \chi(x+\gamma) =-1$  (\cite{berndt,schip3,winter}). 
\QED

\paragraph{Remark 1.} From (\ref{coeff1}) the desired values can be easily read, depending on whether
 $\beta$ is in $\mathcal B_0$ or $\mathcal B_1$ and whether $p^m$ is congruent to $1$ or $3$ modulo $4$ (which determines $\chi_2(-1)$ by the properties of the Jacobi symbol and the quadratic reciprocity law); in particular, we can also read the values for which $i\neq j$, that are not usually explicitly included in the literature; in this case equation (\ref{coeff1}) becomes
\begin{equation}
 \label{mon1}
\frac{1}{4}\left(p^m-2+ \chi_2(-1)\right)~~.
\end{equation}

\begin{theorem}
  \label{lemma3}
The number of representations $R_{p^m}^{(3)}(\beta,i,j)$ of a $\beta\neq 0$ in the field $\mathbb F_{p^m}$
 as the sum of an element in $\mathcal A_j$ and an element in $\mathcal A_i$ is 
\begin{equation}
\label{fieldeven}
 R_{p^m}^{(3)}(\beta,i,j)=\frac{1}{9}(p^m-2-K-\bar K)
\end{equation}
with 
$$
K= \chi_3(\beta)(\omega^{2i}+\omega^{2j})+\omega^{2i+j}-\omega^{2i+2j}\bar\chi_3(\beta)J(\chi_3,\chi_3), 
$$
%

\noindent $J(\chi_3,\chi_3)$ being a Jacobi sum  ($\sum_{c_1+c_2=1} \chi_3(c_1)\chi_3(c_2)$), and
%
depends only on the cubic residuacity of $\beta$.
\end{theorem}

\noindent
{ \sc Proof}. Expanding equation (\ref{count3}), using  $\sum_{x \in \mathbb F_{p^m}} \chi_3(x)=0$ and
 $\sum_{x \in \mathbb F_{p^m}} \chi_3(x) \bar \chi_3(x+\gamma) =-1$, and given that $\chi_3(-1)=1$, we get
$$
\frac{1}{9}\left[p^m-2-\chi_3(\beta)(\omega^{2i}+\omega^{2j})-\bar\chi_3(\beta)(\omega^{i}+\omega^{j})- (\omega^{2i+j}+\omega^{i+2j})
   +\omega^{2i+2j}A(\beta)+\omega^{i+j}\bar A(\beta)\right]
$$
where $A(\beta)=\sum_{x \in  \mathbb F_{p^{m}}}  \chi_3(x) \chi_3(\beta-x)$; 
this summation can be manipulated as 
$$
\sum_{x \in  \mathbb F_{p^{m}}}\chi_3(x) \chi_3(\beta-x)=\chi_3(\beta^2)\sum_{x \in  \mathbb F_{p^{m}}}\chi_3(\beta^{-1}x) \chi_3(1-\beta^{-1}x)=\bar\chi(\beta)A(1)=\bar\chi(\beta)J(\chi_3,\chi_3),
$$
whence the conclusion follows.


%

\paragraph{Remark 2.} The last expression can be further simplified taking into account that 
$J(\chi_3,\chi_3)$ can be computed \cite[Th. 5.21]{lidl} with the Gauss sums of cubic characters (cf. also \cite{schipcubic,schip3}),  as
$$ J(\chi_3,\chi_3) = \frac{G_m^2(1,\chi_3)}{G_m(1,\bar\chi_3)} ~~.  $$ 
In particular, if $p=2$, 
$$\sum_{x \in  \mathbb F_{2^{m}}}  \chi_3(x) \chi_3(x+1) =J(\chi_3,\chi_3)= G_{m}(1,\chi_3) = -(-2)^{m/2}~~.$$ 

\begin{theorem}
If $p\equiv 1$ mod $4$, then $R_{p^m}^{(2)}(0,i,j)$ is $\frac{p-1}{2}$ for $i=j$ or $0$ if $i\neq j$; if $p\equiv 3$ mod $4$, then $R_{p^m}^{(2)}(0,i,j)$ is $\frac{p-1}{2}$ for $i\neq j$ or $0$ if $i=j$.
$R_{p^m}^{(3)}(0,i,j)$ is $\frac{p-1}{3}$ for $i=j$ or $0$ if $i\neq j$.
\end{theorem}
\noindent
{ \sc Proof}. The proof is immediate, taking into account that an element $\alpha$ is in the same coset as $-\alpha$ exactly when $\chi_2(-1)=1$ (resp. $\chi_3(-1)=1$).
\QED
\vspace{5mm}

Working with the characteristic functions, the counterpart of the above theorems would be to multiply $f_{\mathcal B_i}(\mathbf X)$ and $f_{\mathcal B_j}(\mathbf X)$ (or $f_{\mathcal A_i}(\mathbf X)$ and  $f_{\mathcal A_j}(\mathbf X)$) modulo $\mathfrak I_{\mathbf X}$, and then read the coefficients in the output, which involves exactly the same computations as above. But there is more: the characteristic functions 
$f_{\mathcal B_i}(\mathbf X)$ and $f_{\mathcal A_j}(\mathbf X)$ 
satisfy equations of second and third degree, respectively, 
 whose coefficients  are intrinsicly related to the number of representations of the field elements as sums of elements with given quadratic or cubic
 residuacity.

\begin{theorem}
  \label{theo1}
The characteristic functions $f_{\mathcal B_0}(\mathbf X)$ and $f_{\mathcal B_1}(\mathbf X)$ are roots
 of a quadratic equation 
\begin{equation}
  \label{fundeq1}
y^2-\sigma_1 y +\sigma_2=0 \bmod \mathfrak I_{\mathbf X}
\end{equation}
 in the residue ring of multivariate polynomials $\mathbb Z[\mathbf X]/\mathfrak I_{\mathbf X}$.
The sum and the product of the roots (polynomials) are
$$  \left\{ \begin{array}{lcl}
       \sigma_1= -1+\Phi(\mathbf X)  \bmod \mathfrak I_{\mathbf X}\\
       \sigma_2= -\frac{1}{4} \left[p^m \chi_2(-1)-1-\Phi(\mathbf X)(p^m-2+\chi_2(-1)) \right] 
          \bmod \mathfrak I_{\mathbf X}  
\\
      \end{array}  \right. 
$$
 
\end{theorem}

\noindent
{\sc Proof}.
Throughout the proof all the multivariate polynomials should be intended modulo the ideal $\mathfrak I_{\mathbf X}$.
The coefficient $\sigma_1$ is directly obtained from equation (\ref{coeff0}), that is
$$  \sigma_1 = f_{\mathcal B_0}(\mathbf X)+f_{\mathcal B_1}(\mathbf X) = \sum_{\gamma \in \mathbb F_{p^m}^*}
  (I_{\mathcal B_0}(\gamma)+I_{\mathcal B_1}(\gamma)) \mathbf X^{\gamma} = -1+ \Phi(\mathbf X) 
   ~~. $$
The coefficient $\sigma_2$ is computed using the symbols $R^{(2)}_{p^m}(\beta,i,j)$ as follows. 
Starting from
$$ \sigma_2= f_{\mathcal B_0}(\mathbf X)f_{\mathcal B_1}(\mathbf X) = 
 \sum_{\gamma \in \mathbb F_{p^m}^*}  \sum_{\kappa \in \mathbb F_{p^m}^*}
  I_{\mathcal B_0}(\gamma)I_{\mathcal B_1}(\kappa) \mathbf X^{\gamma+\kappa}
~~, $$
and observing that exchanging the summation indices (variables) leaves the result invariant,
we may consider the symmetric summation
$$ \sigma_2 = \frac{1}{2}  \sum_{\gamma \in \mathbb F_{p^m}^*}  \sum_{\kappa \in \mathbb F_{p^m}^*}
\left[I_{\mathcal B_0}(\gamma)I_{\mathcal B_1}(\kappa)+ I_{\mathcal B_0}(\kappa)I_{\mathcal B_1}(\gamma)\right]
  \mathbf X^{\gamma+\kappa} 
 ~~, $$                                         
     and perform the index substitution $\kappa=\beta-\gamma$; the index $\beta$ may be $0$ but it cannot assume the value $\gamma$; thus we may write
$$ \sigma_2 = \frac{1}{2}  \sum_{\gamma \in \mathbb F_{p^m}^*}  
 \sum_{\stackrel{\beta \in \mathbb F_{p^m}}{\beta \neq \gamma}} \left[I_{\mathcal B_0}(\gamma)
  I_{\mathcal B_1}(\beta-\gamma)+ I_{\mathcal B_0}(\beta-\gamma)I_{\mathcal B_1}(\gamma)\right]
  \mathbf X^{\beta} 
~~. $$
Now, in the summation over $\beta$ we separate the term with $\beta=0$ and write
$$ \sigma_2 = \frac{1}{2}  \sum_{\gamma \in \mathbb F_{p^m}^*}  \left\{
 \sum_{\stackrel{\beta \in \mathbb F_{p^m}^*}{\beta \neq \gamma}} \left[I_{\mathcal B_0}(\gamma)
  I_{\mathcal B_1}(\beta-\gamma)+ I_{\mathcal B_0}(\beta-\gamma)I_{\mathcal B_1}(\gamma)\right]
  \mathbf X^{\beta} + C \right\} 
~~. $$
where $C=\left[I_{\mathcal B_0}(\gamma)I_{\mathcal B_1}(-\gamma)+ I_{\mathcal B_0}(-\gamma)I_{\mathcal B_1}(\gamma)\right]= \frac{1-\chi_2(-1)}{2}$, thus exchanging the two summations and, noting that
$ \sum_{\gamma \in \mathbb F_{p^m}^*} C = (p^m-1) C $, we have
$$ \sigma_2 = \frac{p^m-1}{2} C +  \frac{1}{2} \sum_{\beta \in \mathbb F_{p^m}^*} \mathbf X^{\beta} \left\{   
   \sum_{\stackrel{\gamma \in \mathbb F_{p^m}^*}{\gamma \neq \beta}} \left[I_{\mathcal B_0}(\gamma)
  I_{\mathcal B_1}(\beta-\gamma)+ I_{\mathcal B_0}(\beta-\gamma)I_{\mathcal B_1}(\gamma)\right]                   \right\} 
~~. $$
Recalling the definition of $R_{p^m}^{(2)}(\beta,i,j)$, we may write the summation over $\gamma$ 
 as
$$ \sigma_2 = \frac{p^m-1}{2} C +  \frac{1}{2} \sum_{\beta \in \mathbb F_{p^m}^*} \mathbf X^{\beta} \left\{   
   R_{p^m}^{(2)}(\beta,0,1)+R_{p^m}^{(2)}(\beta,1,0) \right\} 
~~. $$
In conclusion, since $R_{p^m}^{(2)}(\beta,0,1)=R_{p^m}^{(2)}(\beta,1,0)$ does not depend on $\beta$, we obtain
$$ \sigma_2 = \frac{p^m-1}{2} \frac{1-\chi_2(-1)}{2} +  R_{p^m}^{(2)}(\beta,0,1) 
 (\Phi(\mathbf X) -1)~~, $$
and, using (\ref{mon1}), we finally obtain
$$  \sigma_2= -\frac{1}{4} \left[p^m \chi_2(-1)-1-\Phi(\mathbf X)(p^m-2+\chi_2(-1)) \right] 
~~. $$
\QED

\begin{theorem}
  \label{theo3}
The characteristic functions $f_{\mathcal A_0}(\mathbf X)$, $f_{\mathcal A_1}(\mathbf X)$, and
 $f_{\mathcal A_2}(\mathbf X)$ are roots of a cubic equation 
\begin{equation}
  \label{fundeq3}
 y^3-\sigma_1 y^2 +\sigma_2 y -\sigma_3 =0 ~~\bmod \mathfrak I_{\mathbf X} ~~,
\end{equation} 
where 
$$   \left\{
    \begin{array}{lcl}
       \sigma_1 &=& \Phi(\mathbf X)-1 \bmod \mathfrak I_{\mathbf X}\\
       \sigma_2 &=&  \frac{1}{3}(p^m-1)(\Phi(\mathbf X)-1) \bmod \mathfrak I_{\mathbf X} 
\\
       \sigma_3 &=&  \frac{1}{27} \left[(\Phi(\mathbf X) -1)^3+(3-3\Phi(\mathbf X)+J(\chi_3,\chi_3)+\bar J(\chi_3,\chi_3))(p^m-\Phi(\mathbf X)) \right] \bmod \mathfrak I_{\mathbf X} 
\\
    \end{array}  \right.
$$    
\end{theorem}

\noindent
{\sc Proof}. Throughout the proof all the multivariate polynomials should be intended modulo the ideal $\mathfrak I_{\mathbf X}$.
The coefficient $\sigma_1$ is easily computed as
$$ f_{\mathcal A_0}(\mathbf X)+ f_{\mathcal A_1}(\mathbf X)+f_{\mathcal A_2}(\mathbf X) = \sum_{\gamma \in \mathbb F_{p^m}^*}
  (I_{\mathcal A_0}(\gamma)+ I_{\mathcal A_1}(\gamma)+I_{\mathcal A_2}(\gamma)) \mathbf X^{\gamma} = -1+ \Phi(\mathbf X) ~~, $$
because $I_{\mathcal A_0}(\gamma)+I_{\mathcal A_1}(\gamma)+I_{\mathcal A_2}(\gamma)=1$ and equation (\ref{coeff0}) is used.
The elementary symmetric function $\sigma_2$ is the sum
$$ f_{\mathcal A_0}(\mathbf X)f_{\mathcal A_1}(\mathbf X)+ 
   f_{\mathcal A_1}(\mathbf X)f_{\mathcal A_2}(\mathbf X)+
   f_{\mathcal A_2}(\mathbf X)f_{\mathcal A_0}(\mathbf X)  ~~, 
$$
then 
 we need to compute the summation
$$ \sigma_2 =  \sum_{\gamma \in \mathbb F_{p^m}^*}
 \sum_{\eta \in \mathbb F_{p^m}^*}
  (I_{\mathcal A_0}(\gamma)I_{\mathcal A_1}(\theta)+ I_{\mathcal A_1}(\gamma)I_{\mathcal A_2}(\theta)+
  I_{\mathcal A_2}(\gamma)I_{\mathcal A_0}(\theta)) \mathbf X^{\gamma+\theta} 
~~. $$
Using (\ref{fieldeven}), 
 we get
$$ \sigma_2=  - \frac{1}{3}(p^m-1) +\frac{1}{3}(p^m-1)  \Phi(\mathbf X) 
~~. $$  
Lastly, the elementary symmetric function 
 $\sigma_3=f_{\mathcal A_0}(\mathbf X)f_{\mathcal A_1}(\mathbf X) f_{\mathcal A_2}(\mathbf X)$ is given by
  the summation   
$$ \sigma_3 =  \sum_{\gamma \in \mathbb F_{p^m}^*}
 \sum_{\theta \in \mathbb F_{p^m}^*} \sum_{\kappa \in \mathbb F_{p^m}^*}
  I_{\mathcal A_0}(\gamma)I_{\mathcal A_1}(\theta) I_{\mathcal A_2}(\kappa)
   \mathbf X^{\gamma+\theta+\kappa}
 ~~. $$
Expanding the product of the indicator functions and performing the summations, most of the $27$ sums
 are cancelled, and it remains to compute the following:
$$  \frac{1}{27}  \sum_{\gamma \in \mathbb F_{q}^*}
 \sum_{\theta \in \mathbb F_{p^m}^*} \sum_{\kappa \in \mathbb F_{p^m}^*} \left(1-3 \chi_3(\gamma) \bar \chi_3(\theta)
 +\chi_3(\gamma) \chi_3(\theta) \chi_3(\kappa)  + \bar \chi_3(\gamma) \bar \chi_3(\theta) 
  \bar \chi_3(\kappa)  \right) \mathbf X^{\gamma+\theta+\kappa}  ~~.  $$
To complete the task we then need to compute only three summations.
\begin{enumerate}
 \item The summation
$$ \sum_{\gamma \in \mathbb F_{p^m}^*}
 \sum_{\theta \in \mathbb F_{p^m}^*} \sum_{\kappa \in \mathbb F_{p^m}^*} \mathbf X^{\gamma+\theta+\kappa}  = ( \Phi(\mathbf X)-1)^3 ~~; $$
    is easily obtained, because the three summations on $\gamma$, $\theta$, and $\kappa$ can be performed independently.
  \item The summation  $\sty \sum_{\gamma \in \mathbb F_{p^m}^*}
 \sum_{\theta \in \mathbb F_{p^m}^*} \sum_{\kappa \in \mathbb F_{p^m}^*} \chi_3(\gamma) \bar \chi_3(\theta)
 \mathbf X^{\gamma+\theta+\kappa}$ is computed by extending the sum range to include $0$; this is  done using the function $\delta(\kappa)$ which is $1$ if $\kappa=0$, and is $0$
 otherwise, thus the summation is $\sty \sum_{\gamma \in \mathbb F_{p^m}}
 \sum_{\theta \in \mathbb F_{p^m}} \sum_{\kappa \in \mathbb F_{p^m}} \chi_3(\gamma) \bar \chi_3(\theta)
 (1-\delta(\kappa)) \mathbf X^{\gamma+\theta+\kappa} $ which splits into two summations
$$  \sum_{\gamma \in \mathbb F_{p^m}}
 \sum_{\theta \in \mathbb F_{p^m}} \sum_{\kappa \in \mathbb F_{p^m}} \chi_3(\gamma) \bar \chi_3(\theta)
\mathbf X^{\gamma+\theta+\kappa} - \sum_{\gamma \in \mathbb F_{p^m}}
 \sum_{\theta \in \mathbb F_{p^m}} \sum_{\kappa \in \mathbb F_{p^m}} \chi_3(\gamma) \bar \chi_3(\theta)
 \delta(\kappa)\mathbf X^{\gamma+\theta+\kappa} ~~. $$
In the triple summations, the sum over $\kappa$ can be performed independently and gives
 $\Phi(\mathbf X)$ for the first, and simply $1$ for the second. Thus we have
$$ (\Phi(\mathbf X)-1)  \sum_{\gamma \in \mathbb F_{p^m}}
 \sum_{\theta \in \mathbb F_{p^m}}  \chi_3(\gamma) \bar \chi_3(\theta)
\mathbf X^{\gamma+\theta} ~~.  $$
The double summation can be easily evaluated with the substituttion $\theta=\beta-\gamma$
$$  \sum_{\beta \in \mathbb F_{p^m}}  \sum_{\gamma \in \mathbb F_{p^m}} \chi_3(\gamma) \bar \chi_3(\beta-\gamma) \mathbf X^{\beta}  = p^m-1 - \sum_{\beta \in \mathbb F_{p^m}^*}   \mathbf X^{\beta} = p^m-\Phi(\mathbf X)  ~~, $$
because the sum over $\gamma$ assumes only two values, namely $-1$ if $\beta \neq 0$ and $p^m-1$ if $\beta=0$.
In conclusion we obtain 
$$(\Phi(\mathbf X)-1)(p^m-\Phi(\mathbf X)) ~~. $$ 
\item The sums in the triple summation
$$ \sum_{\gamma \in \mathbb F_{p^m}^*}
 \sum_{\theta \in \mathbb F_{p^m}^*} \sum_{\kappa \in \mathbb F_{p^m}^*} \chi_3(\gamma)  \chi_3(\theta)  \chi_3(\kappa)
 \mathbf X^{\gamma+\theta+\kappa} = \sum_{\beta \in \mathbb F_{p^m}}
 \sum_{\theta \in \mathbb F_{p^m}} \sum_{\gamma \in \mathbb F_{p^m}} \chi_3(\gamma)  \chi_3(\theta)  \chi_3(\beta-(\gamma+\theta))  \mathbf X^{\beta}  ~~, $$ 
have been extended throughout $\mathbb F_{p^m}$ as $\chi_3(0)=0$, together with the substitution
 $\kappa=\beta-(\gamma+\theta)$. Now, the summation over $\gamma$ has two values, namely
 $0$ if $\beta=\theta$, and $\bar \chi_3(\beta-\theta) A(1)$ if $\beta \neq \theta$, $A(1)$ being as above $\sum_{x \in  \mathbb F_{p^{m}}}  \chi_3(x) \chi_3(1-x)=J(\chi_3,\chi_3)$,
 therefore we obtain 
$$\sum_{\beta \in \mathbb F_{p^m}}
 \sum_{\stackrel{\theta \in \mathbb F_{p^m}}{\theta \neq \beta}} \chi_3(\theta) \bar \chi_3(\beta-\theta) A(1)
 \mathbf X^{\beta} = \sum_{\beta \in \mathbb F_{p^m}}
 \sum_{\theta \in \mathbb F_{p^m}} \chi_3(\theta) \bar \chi_3(\beta-\theta)A(1)
 \mathbf X^{\beta} = A(1) (p^m-\Phi(\mathbf X))  $$ 
since the restriction $\theta \neq \beta$ can be removed and the summation over $\theta$ is $-1$ if $\beta\neq 0$ or $p^m-1$ if $\beta= 0$. 

\end{enumerate}     

In conclusion, collecting the results we obtain 
$$  \sigma_3 = \frac{1}{27} \left((\Phi(\mathbf X) -1)^3-3(\Phi(\mathbf X)-1)(p^m-\Phi(\mathbf X)) +
         [J(\chi_3,\chi_3)+\bar J(\chi_3,\chi_3)] (p^m-\Phi(\mathbf X))\right)  
  $$   
\QED

\paragraph{Remark 3.} Note that $J(\chi_3,\chi_3)+\bar J(\chi_3,\chi_3)$ is always an integer, being twice the sum of real parts of cubic roots of unity. 


\paragraph{Remark 4.} Even though its derivation involved handling products of three characters, the expression of $\sigma_3$ only involves the Jacobi sum $A(1)$, i.e. fundamentally only the number of representations as the sum of two elements of two given cosets. 

\section{Connections with other problems}

In this section we point out a sort of duality relationship with the following problem.

Suppose that we have $t$ elements of a finite field $\mathbb F_{p^m}$ all belonging to one of the cosets
 determined by the character partition. We would like to know how many $\beta$s there are in the field such
 that, adding $\beta$ to all the $t$ elements, we get $t$ elements still belonging to a common coset.
If the character has order $n$, we let $N_{p^m}^{(n)}(t)$ be the number of $\beta$s; i.e. it is the
 number of solutions $\beta$ of a system of $t$ equations in $\mathbb F_{p^m}$ of the form
\begin{equation}
   \label{sys1}
  \left\{ \begin{array}{l}
        \alpha^j z_1^n +\beta= \alpha^k y_1^n \\
       \alpha^j  z_2^n +\beta= \alpha^k y_2^n \\
 ~~~~~~\vdots    \\
        \alpha^j z_t^n +\beta= \alpha^k y_t^n \\ 
    \end{array}  \right.
\end{equation}
where $\alpha^j z_1^n, \alpha^j  z_2^n, \cdots, \alpha^j z_t^n$ are given and distinct, $\alpha$ being a primitive element, whereas the elements $y_i$s must be chosen in the field to satisfy the system, and the $n$ values $\{0,1,\ldots,n-1\}$ for $k$ and $j$ are all considered.
 However, we may assume $j=0$, since dividing each equation by $\alpha^j$, and setting
  $\beta'=\beta \alpha^{-j}$ and $k'=k-j \bmod n$, we see that the number of solutions 
  of the system is independent of $j$.

An explicit solution when the character is quadratic or cubic can be obtained, again by means of the indicator functions. For example, if we have a cubic character over $\mathbb F_{2^m}$, given a $z_i$ we can partition the elements $\beta\neq z_i^3$ in $\mathbb F_{2^m}$      
 into subsets depending on the $k\in\{0,1,2\}$ such that $\chi(\beta+z_i^3)=\omega^k$. 
Therefore, a solution of (\ref{sys1}) for a fixed $k$ and $j=0$ is singled out by the
 product
$$  \prod_{i=1}^t I_{\mathcal A_k}(\beta+z_i^3)=\frac{1}{3^t} [1+\sum_{i=1}^t \sigma_i^{(k)} ]   ~~, $$
where each $\sigma_i^{(k)}$ is a homogeneous sum of monomials which are products of $i$ characters of the form $\chi(\beta+z_h^3)$ or $\bar\chi(\beta+z_h^3)$.
Thus $N_{2^m}^{(3)}(t)$ is
\begin{equation}
   \label{main1}
   N_{2^m}^{(3)}(t)= \sum_{\stackrel{\beta \in \mathbb F_{2^m}}{\beta\not \in\{z_i^3\}}} 
        \left[ \prod_{i=1}^t I_{\mathcal A_0}(\beta+z_i^3) + \prod_{i=1}^t I_{\mathcal A_1}(\beta+z_i^3)+\prod_{i=1}^t I_{\mathcal A_2}(\beta+z_i^3) \right] ~~.
\end{equation}
The $z_i$ are excluded from the sum, since $z_i^3+z_i^3=0$ does not belong to any coset.

In \cite{schip7}, which deals with this problem to analyse the success rate of the Cantor-Zassenhaus
 polynomial factorization algorithm, we computed exactly some of the above expressions for small values
 of $t$, and gave bounds for more general cases. In particular, we found that
$$   1+\max_{z_1\neq z_2\neq z_3} N_{2^m}^{(3)}(3)= \left\{
     \begin{array}{l}
       \frac{1}{9}(2^m+2^{m/2}-2) ~~~~\mbox{for}~~ m/2 \mbox{ even}   \\
       \frac{1}{9}(2^m+2^{m/2+1}+1) ~~~~\mbox{for}~~ m/2 \mbox{ odd}
     \end{array} \right.
        ~~ $$   
and
$$
1+\max_{z_1\neq z_2\neq z_3} N_{p^m}^{(2)}(3)=\left\{\begin{array}{l}
\frac{1}{4}(p^m-1)~~~~\mbox{for}~~p=4k+1\\
\frac{1}{4}(p^m+1)~~~~\mbox{for}~~p=4k+3,~~m\ \mbox{odd}\\
\frac{1}{4}(p^m-1)~~~~\mbox{for}~~p=4k+3,~~m\ \mbox{even}
\end{array}\right.
.$$

Recalling that $R_{p^m}^{(n)}(\beta,i,j)$, $n=2,3$, denotes 
the number of representations of a $\beta\neq 0$ in a finite field $\mathbb F_{p^m}$ 
 as the sum of two element belonging to two cosets indexed by $i$ and $j$ in the partition given by a character of order $n$, 
we find the remarkable identities:
$$
\max_{i,j,\beta} R_{2^m}^{(3)}(\beta,i,j)=1+\max_{z_1\neq z_2\neq z_3} N_{2^m}^{(3)}(3)
$$
and
$$
\max_{i,j,\beta} R_{p^m}^{(2)}(\beta,i,j)=1+\max_{z_1\neq z_2\neq z_3} N_{p^m}^{(2)}(3)~~.
$$

In \cite{schip7} this quasi-duality had the following interesting interpretation: the maximum $t$ such that it is still possible to fail to split a polynomial of degree $t$ with two attempts is equal to the maximum number of attempts to split a polynomial of degree $3$. 

\section{Acknowledgments}
%
The Research was supported in part by the Swiss National Science
Foundation under grant No. 132256.


\begin{thebibliography}{99}

\bibitem{berndt}
    B. Berndt, R.J. Evans, H. Williams, 
    {\em Gauss and Jacobi Sums},
   Wiley, 1998.
\bibitem{schip7}
     M. Elia, D. Schipani, 
     Improvements on the Cantor-Zassenhaus Factorization Algorithm,
      {\it www.arxiv.org}, 2011.
\bibitem{schipcubic}
     M. Elia, D. Schipani, 
     Gauss sums of cubic characters over $GF(p^r)$, $p$ odd,
      {\it www.arxiv.org}, 2011.

\bibitem{lidl}
    R. Lidl, H. Niederreiter,
    {\em  Finite Fields},
    Cambridge Univ. Press, 1997.
\bibitem{elia2}
     C. Monico, M. Elia,
       {Note on an Additive Characterization of Quadratic Residues Modulo $p$,}
       {\em Journal of Combinatorics, Information \& System Sciences}, vol. 31, 2006, pp.209-215.
\bibitem{monico}
    C. Monico, M. Elia,
     An Additive Characterization of Fibers of Characters on $\mathbb F_p^*$,  
     {\em International Journal of Algebra}, Vol. 1-4, n.3, 2010, pp.109-117.
\bibitem{perron}
    O. Perron,
     Bemerkungen uber die Verteilung der quadratischen Reste,  
     {\em Mathematische Zeitschrift}, Vol. 56, 1952, pp.122-130.
\bibitem{raymond}
     D. Raymond,
     {\em An Additive Characterization of Quadratic Residues},
      Master Degree thesis, Texas Tech University (Lubbock), 2009.
\bibitem{schip3}
    D. Schipani, M. Elia,
     Gauss Sums of the Cubic Character over $\mathbb F_{2^m}$: an elementary derivation,
     \emph{Bull. Polish Acad. Sci.Math.}, 59, 2011, pp.11-18.
\bibitem{winter}
     A. Winterhof,
     On the Distribution of Powers in Finite Fields,
   {\em Finite Fields and Their applications}, 4, 1998, pp.43-54.
\end{thebibliography}
\end{document}